\documentclass[11pt]{article}
\usepackage{graphics}
\usepackage[demo]{graphicx}
\usepackage{epsfig}
\usepackage{pstricks}
\usepackage[normal]{subfigure}
\usepackage[latin1]{inputenc}
\usepackage[english,francais]{babel}
\usepackage{relsize,exscale}
\usepackage{makeidx}
\usepackage{enumitem}
\usepackage{amsfonts,amssymb,amsmath}
\usepackage{graphicx}
\usepackage{color}
\usepackage{multirow}
\usepackage{mathrsfs}
\usepackage[normalem]{ulem}
\usepackage{cancel}
\usepackage{bbm}
\usepackage{empheq}
\newenvironment{prooff}{{\it Proof :}}{\hfill\rule{2mm}{2mm}\vskip3mm\par}
\newtheorem{theorem}{Theorem}[section]
\newtheorem{lemma}[theorem]{Lemma}
\newtheorem{proposition}[theorem]{Proposition}
\newtheorem{corollary}[theorem]{Corollary}
\newtheorem{e-definition}[theorem]{Definition\rm}
\newtheorem{remark}{\it Remark\/}

\usepackage{color}
\definecolor{dred}{rgb}{0.92,0,0}
\definecolor{dgreen}{rgb}{0,0.92,0}
\definecolor{dblue}{rgb}{0,0,0.92}
\definecolor{dyellow}{rgb}{0.95,0.95,0}

\newcommand{\R}{\mathbb{R}}
\def\h{{\sf h}}
\def\PiS{\pi_{\mathcal{S}}(v)}
\def\PiSP{\pi_{\mathcal{S}}(v)(P)}
\newcommand{\N}{\mathbb{N}}
\def\D{\displaystyle}
\newcommand{\hs}{\hspace{0.1cm}}

\newcommand{\sa}{\\ [0.2cm]}
\def\restriction#1#2{\mathchoice
              {\setbox1\hbox{${\displaystyle #1}_{\scriptstyle #2}$}
              \restrictionaux{#1}{#2}}
              {\setbox1\hbox{${\textstyle #1}_{\scriptstyle #2}$}
              \restrictionaux{#1}{#2}}
              {\setbox1\hbox{${\scriptstyle #1}_{\scriptscriptstyle #2}$}
              \restrictionaux{#1}{#2}}
              {\setbox1\hbox{${\scriptscriptstyle #1}_{\scriptscriptstyle #2}$}
              \restrictionaux{#1}{#2}}}
\def\restrictionaux#1#2{{#1\,\smash{\vrule height .8\ht1 depth .85\dp1}}_{\,#2}}
\usepackage{multirow}
\graphicspath{
{./Figures/}
{Figures/}
{./}
}
\title{A refined first-order expansion formula in $\R^n$:\\
Application to interpolation and finite element error estimates
}
\author{Jo\"el Chaskalovic \thanks{D'Alembert,
Sorbonne University, Paris, France, (Email: \emph{(corresp.)} jch1826@gmail.com)}
\qquad
Franck Assous
\thanks{
Department of Mathematics, Ariel University, Israel, (Email: franckassous55@gmail.com)}
}
\date{}
\begin{document}
\maketitle
\selectlanguage{english}
\begin{abstract}
\noindent The aim of this paper is to derive a refined first-order expansion formula in $\R^n$, the goal being to get an optimal reduced remainder, compared to the one obtained by usual Taylor's formula. For a given function, the formula we derived is obtained by introducing a linear combination of the first derivatives, computed at $n+1$ equally spaced points. We show how this formula can be applied to two important applications: the interpolation error and the finite elements error estimates. In both cases, we illustrate under which conditions a significant improvement of the errors can be obtained, namely how the use of the refined expansion can reduce the upper bound of error estimates. 
\end{abstract}
\noindent {\em keywords}: Taylor's theorem, interpolation error estimates, finite element, approximation error estimates.
\section{Introduction}\label{intro}
\noindent Evaluating and improving the accuracy of approximation are very difficult problems in numerical analysis. In this article, we are interested in describing a new point of view on approaching the topic. In particular, we are concerned with the difficulty of accurately determining the error estimate of numerical methods applied to partial differential equations.\\

\noindent From a mathematical point of view, the origin of this problem \cite{AsCh2014} can be found in Rolle's theorem, and therefore, in Lagrange and Taylor's theorems \cite{Atkinson}: this comes from the existence of a non unique unknown point which appears in the remainder of Taylor's expansion, as the heritage of Rolle's theorem, leading to a kind of ``uncertainty''.\\

\noindent Now, let us consider more specifically finite element problems. Basically because of this uncertainty, most of the results focus on the asymptotic behavior of the error estimates which strongly depends on the interpolation error (see for example \cite{ChaskaPDE}, \cite{ArXiv_JCH}). Indeed, error estimates generally consider, for a given norm, the asymptotic behavior of the difference between the exact and the approximate solution, as the mesh size $h$ tends to zero (cf. \cite{RaTho82}).\\

\noindent  However, several approaches have been proposed to investigate how to improve the accuracy of approximation. For example, in
the framework of numerical integration, we refer the reader to \cite{Barnett_Dragomir}, \cite{Cerone} or \cite{Dragomir_Sofo}, and references therein. From another point of view, due to the lack of information, heuristic methods were considered, basically based on a probabilistic approach, see for instance \cite{Abdulle}, \cite{Hennig}, \cite{Oates} or \cite{ChAs20} and \cite{CMAM2}. This allows to compare different numerical methods, and more precisely finite element, for a given fixed mesh size, see \cite{MMA2021}-\cite{AxiChAs21}.\\

\noindent Nevertheless, Taylor's formula introduces an unknown point preventing to precisely determine the interpolation error, and consequently the approximation error of a given numerical method. Hence, the question still remains whether the error upper bounds are as small as possible. Therefore, we focus in this article on the values of the numerical constants which appear in such estimates, aiming to reduce them as much as possible.\\

\noindent To this end, we propose to refine the first order Taylor approximation formula of a function $f$ in $\R^n$, at a given point $a$, using $f(a)$ and its derivative. To do this, we consider more known values of the first order derivative on intermediary equidistant points between $a$ and another given point $b$. These values of the first order derivative are conveniently weighted in order to diminish the associated reminder. Then, we study the resulting properties in the interpolation error estimates and in Lagrange finite element error estimates.\\

\noindent The paper is organized as follows. In Section \ref{B}, we present the main result of this paper which deals with a new refined first-order
expansion formula in $\R^n$. Section \ref{applInter} investigates the consequences in interpolation errors. The case of dimension one and dimension $n>1$ are separately investigated.  Application to finite elements errors estimates is studied in section \ref{appliFE}. Some remarks are proposed. Concluding remarks follow.

\section{A new first order expansion formula in $\R^n$}\label{B}
\noindent We consider a non-empty bounded and simply connected open domain $U\subset\R^n, (n\in\N^{*})$, and  a given point $a\in U$. We also consider a function $f$ which is twice differentiable on $U$. \sa
We recall that the second order differential $D^2\!f(x)$ at a given point $x\in U$ belongs to the vector space of linear and continuous forms $\mathcal{L}(\R^n,\mathcal{L}(\R^n,\R))$
which can be identified to the space of bilinear, symmetric and continuous form $\mathcal{B}(\R^n\times\R^n,\R)$.\sa
In other words, we have
$$
\forall \,(h,h')\in \R^n\times\R^n: \hs D^{2}\!f(x).(h,h') = \sum_{i,j=1}^{n}h_ih'_j\frac{\partial^2f}{\partial x_i\partial x_j}(x)\,.
$$
Above, and in the rest of this article, we adopt the following writing convention: For any linear form $L$, we denote by $L.(h)$ or $L(h)$ the action of $L$ on a given vector $h\in\R^n$, and by
$b.(h,h')$ or $b(h,h')$ the action of a bilinear form $b$ on a given couple of vectors $(h,h')\in\R^n\times\R^n$ (see for example A. Avez \cite{Avez}).\sa
Moreover, we assume that there exist two constants $(m_2,M_2)\in\R^2$ such that
\begin{equation}\label{D2f_bounded}
\forall \, x\in U, \forall (h,h')\in \R^n\times\R^n: m_2\|h\|\|h'\| \hs \leq \hs D^2\!f(x).(h,h') \leq M_2\|h\|\|h'\|,
\end{equation}
where $\|.\|$ denotes a given norm defined on $\R^n$. Remark that condition (\ref{D2f_bounded}) can be interpreted by using the natural norm $|||.|||$ of a bilinear and continuous form $b$, defined by:
\begin{equation}\label{norme_forme_bili}
\D \forall \, b \in \mathcal{B}(\R^n\times\R^n,\R): |||b||| \equiv \sup_{(h,h')\in \R^{n^*}\!\!\times\R^{n^*}}\frac{|b(h,h')|}{\|h\|.\|h'\|}.
\end{equation}
Then, applying this definition to the bilinear form $b=D^2\!f(x)$,  inequality (\ref{D2f_bounded}) can be written as:
$$
\exists \mathcal{M}_2 > 0, \forall \, x \in U: |||D^2\!f(x)||| \leq \mathcal{M}_2, (\mathcal{M}_2=\max(|m_2|,|M_2|)).
$$
We first begin with a straight consequence of the classical first-order Taylor formula. 
\begin{proposition}\label{Prop1}
Let $f$ be a twice-differentiable function defined on a non-empty bounded  and simply connected open set $U \in \R^n$. Then, we have
\begin{equation}\label{Taylor_1_V1}
\D f(a+h) = f(a) + D\!f(a).(h) + \|h\|\epsilon_{a,1}(h),
\end{equation}
where $\|.\|$ denotes a given norm on $\R^n$, the remainder $\epsilon_{a,1}(h)$ satisfiying:
\begin{equation}\label{epsilon1_bounded}
\frac{\|h\|}{2}m_2 \leq \epsilon_{a,1}(h) \leq \frac{\|h\|}{2}M_2.
\end{equation}
\end{proposition}
\begin{prooff}
For a given twice-differentiable function $f$, we have the classical Taylor's expansion (see for example \cite{Taylor_Rn}):
$$
\D\exists\,\xi_{a,h}\in \, ]a,a+h[ \hs \mbox{ such that} \hs f(a+h) = f(a) + D\!f(a).(h) + \frac{1}{2}\,D^2\!f(\xi_{a,h}).(h,h),
$$
where $]a,a+h[$ denotes the open line segment bounded by the two points $a$ and $a+h$. For $h$ sufficiently small, $]a,a+h[ \in U$.
Then, by comparing this expression with (\ref{Taylor_1_V1}), we have
\begin{equation}\label{epsilon_D2f}
\|h\|\epsilon_{a,1}(h) = \frac{1}{2}\,D^2\!f(\xi_{a,h}).(h,h).
\end{equation}
So, using (\ref{D2f_bounded}) with $h'=h$, we get (\ref{epsilon1_bounded}) from (\ref{epsilon_D2f}).
\end{prooff}
Now, to derive the main result, let us first introduce the function $\phi$ defined as follows:
\begin{equation}\label{phi}
\begin{array}{r c c c l}
    \phi & : & [0,1] & \longrightarrow & \mathbb{R} \\
         &   &     t & \longmapsto & D\!f(a + th).(h)\,.
\end{array}
\end{equation}
First, we remark that $\phi(0) = D\!f(a).(h)$ and that $\phi(1) = D\!f(a+h).(h)$.
Moreover, the remainder $\epsilon_{a,1}(h)$ introduced in (\ref{Taylor_1_V1}) satisfies the following result:
\begin{proposition}\label{Prop2}
Let $\epsilon_{a,1}(h)$ defined by the first order Taylor expansion (\ref{Taylor_1_V1}). Then, we have
\begin{equation}\label{Epsilon1}
\|h\|\epsilon_{a,1}(h) = \int_{0}^{1}(1-t)\phi'(t)\,dt.
\end{equation}
\end{proposition}
\begin{prooff}
Using the first-order Taylor's formula with the integral form of the remainder gives:
\begin{equation}\label{Taylor_Integral_V0}
f(a+h) = f(a) + D\!f(a).(h) + \int_{0}^{1}(1-t)D^2\!f(a+th).(h,h)\,dt.
\end{equation}
Now, deriving function $\phi$ defined by (\ref{phi}) with respect to $t$, we obtain:
\begin{eqnarray*}
\phi'(t) & = & \frac{d}{dt}\biggl[D\!f(a + th).(h)\biggr] = \frac{d}{dt}\sum_{i=1}^{n}\frac{\partial f}{\partial x_i}(a + th)h_i, \\[0.2cm]
 & = & \frac{d}{dt}\sum_{i=1}^{n}\frac{\partial f}{\partial x_i}(a_1 + th_1,\dots, a_n + th_n)h_i.
\end{eqnarray*}
Using standard rule of derivation \cite{Taylor_Rn} to derive a function of $n$ variables, each depending on $t$, we get:
\begin{eqnarray*}
\phi'(t) & = & \sum_{j=1}^{n}\biggl[\sum_{i=1}^{n}\frac{\partial^2 f}{\partial x_j\partial x_i}(a_1 + th_1,\dots, a_n + th_n)h_ih_j\biggr], \\[0.2cm]
& = & D^{2}\!f(a + th).(h,h),
\end{eqnarray*}
and finally, the first order Taylor's formula (\ref{Taylor_Integral_V0}) leads to:
$$
f(a+h) = f(a) + D\!f(a).(h) + \int_{0}^{1}(1-t)\phi'(t) dt.
$$
As consequence, the remainder $\epsilon_{a,1}(h)$ in (\ref{Taylor_1_V1}) satisfies (\ref{Epsilon1}).
\end{prooff}
Consider now a given $m \in \mathbb{N}^{*}$, for which we define $\epsilon_{a,m+1}(h)$ by 
\begin{equation}\label{Def_epsilon_m+1}
f(a+h) = f(a) + \sum \limits_{k=0}^{m} \omega_{k}(m)D\!f\!\left(a + \frac{kh}{m}\right)\!\!.(h) + \|h\|\epsilon_{a,m+1}(h),
\end{equation}
where  $\omega_{k}(m)$, $0 \leq k \leq m$, denote real weights we want  to determine to get a corresponding remainder $\epsilon_{a,m+1}(h)$
as small as possible. Since we assumed that the domain $U$ is open and simply connected, the set of uniformly distributed points $(x_k)_{k=1,m}$ defined by $\D x_k=a + \frac{kh}{m}$ belongs to the segment $[a,a+h]$.\\

\noindent The following theorem constitues our main result, and is devoted to the refined first-order expansion in $\R^n$.
\begin{theorem}\label{theorem_1}
Let $f$ be a twice differentiable real mapping defined on a non-empty bounded and simply connected open $U\subset\R^n$, such that (\ref{D2f_bounded}) holds, and let $a\in U$ be a given point. \sa
If the weights $\omega_k(m), (k=0,m),$ satisfy: $$\D \sum \limits_{k=0}^{m} \omega_{k}(m) = 1,$$
then, the following refined first-order expansion formula holds:
\begin{equation}\label{Generalized_Taylor}
f(a+h) = f(a) + \left(\frac{D\!f(a) + D\!f(a+h)}{2m} + \frac{1}{m}\sum \limits_{k=1}^{m-1} D\!f\!\left(a + \frac{kh}{m}\right)\!\!\right)\!.(h) + \|h\|\epsilon_{a,m+1}(h),
\end{equation}
where the remainder $\epsilon_{a,m+1}(h)$ satisfies
\begin{equation}\label{majoration_epsilon_m+1}
\D \frac{\|h\|}{8m}(m_2-M_2)\leq \epsilon_{a,m+1}(h) \leq \frac{\|h\|}{8m}(M_2-m_2).
\end{equation}
Moreover, for a uniform distribution of points, this result is optimal in the sense that the weights $\omega_k(m)$ in (\ref{Generalized_Taylor}) guarantee the remainder $\epsilon_{a, m+1}(h)$ to be minimal.
\end{theorem}
\newpage
\begin{remark}
$\frac{}{}$
\begin{enumerate}
\item In the theorem above, we use the convention that, when $m=1$, the sum $\displaystyle  \sum \limits_{k=1}^{m-1}$ is identically zero.
\item Remark that the error bounds in estimate (\ref{majoration_epsilon_m+1}) are $2m$ lower that the usual ones involved in estimate (\ref{epsilon1_bounded}).
\end{enumerate}
\end{remark}
In order to prove the theorem \ref{theorem_1}, we will need the following lemma obtained in  \cite{ChaskJams}:
\begin{lemma} \label{formula_1}
Let $u$ be a continuous function on $\R$, and let $(a_{k})_{k=0,m}, (m\in\N^*)$ be a sequence of real numbers. We have the following formula:
$$
 \sum \limits_{k=0}^{m-1} \int_{k}^{m}a_{k} u(t)\,dt = \sum \limits_{k=0}^{m-1} \int_{k}^{k+1}S_{k} u(t)\,dt,
$$
where
$$
\D S_{k} = \sum \limits_{j=0}^{k} a_{j}.
$$
\end{lemma}
Let us now prove Theorem \ref{theorem_1}.\sa
\begin{prooff} 
From (\ref{Taylor_1_V1}) and (\ref{Def_epsilon_m+1}), we have
\begin{equation}\label{FF1}
\D D\!f(a).(h) + \|h\|\epsilon_{a,1}(h) = \sum \limits_{k=0}^{m} \omega_{k}(m)D\!f\!\biggl(\!a+\frac{kh}{m}\biggr)\!.(h) + \|h\|\epsilon_{a,m+1}(h),
\end{equation}
that can be written, using the function $\phi$ introduced in (\ref{phi}),
\begin{equation}\label{FF2}
\phi(0) + \|h\|\epsilon_{a,1}(h) = \sum \limits_{k=0}^{m} \omega_{k}(m)\phi\biggl(\frac{k}{m}\biggr) + \|h\|\epsilon_{a,m+1}(h).
\end{equation}
Using now Proposition \ref{Prop2}, we get from (\ref{FF2}) that the remainder $\epsilon_{a,m+1}(h)$ satisfies:
$$
\|h\|\epsilon_{a,m+1}(h) =  \phi(0) + \int_{0}^{1}(1-t)\phi'(t)\,dt - \sum \limits_{k=0}^{m} \omega_{k}(m)\phi\biggl(\frac{k}{m}\biggr), 
$$
that can be written as (see details in \cite{ChaskJams}, pages 5 and 6, Eq. (17)-(19))
\begin{equation}
\|h\|\epsilon_{a,m+1}(h) =  \biggl(1-\sum_{k=0}^{m}\omega_k(m)\biggr)\phi(1) - \sum_{k=0}^{m-1}\int_{\frac{k}{m}}^{\frac{k+1}{m}}t\phi'(t)\, dt + \sum_{k=0}^{m-1}\int_{\frac{k}{m}}^{1}\omega_k(m)\phi'(t)\, dt. \label{FF3.2}
\end{equation}
Let us now transform the last integral of (\ref{FF3.2}) by the help of Lemma \ref{formula_1}. We set:
$$\D \forall\, t\in\R: u(t) = \phi'\biggl(\frac{t}{m}\biggr), \hs a_{k} = \omega_{k}(m), \mbox{ and } \hs S_{k} = \sum \limits_{j = 0}^{k} \omega_{j}(m) = S_{k}(m).$$
Then, after substitution, (\ref{FF3.2}) becomes:
\begin{equation}\label{FF4}
\|h\|\epsilon_{a,m+1}(h) = \biggl(1-\sum_{k=0}^{m}\omega_k(m)\biggr)\phi(1) + \sum \limits_{k=0}^{m-1} \int_{\frac{k}{m}}^{\frac{k+1}{m}}(S_{k}(m)-t)\phi'(t)\, dt.
\end{equation}
From now on, let us assume for simplicity (see Remark \ref{Relation_Poids} below) that:
\begin{equation}\label{poids_norme}
\D \sum \limits_{k=0}^{m} \omega_{k}(m) = 1,
\end{equation}
so (\ref{FF4}) can be written
\begin{equation}\label{FF5}
\|h\|\epsilon_{a,m+1}(h) = \sum \limits_{k=0}^{m-1} \int_{\frac{k}{m}}^{\frac{k+1}{m}}(S_{k}(m)-t)\phi'(t)\, dt.
\end{equation}
Using estimate (\ref{D2f_bounded}) of the second order differential of $f$, and the definition (\ref{phi}) of $\phi(t)$, we get:
\begin{equation}\label{II.1}
m_2\|h\|^2 \leq \phi'(t) = D^2\!f(a+th).(h,h) \leq M_2\|h\|^2.
\end{equation}
So, to derive a double inequality on $\epsilon_{a,m+1}(h)$, we split the integral in (\ref{FF5})
\begin{equation}\label{FF6}
\int_{\frac{k}{m}}^{\frac{k+1}{m}}(S_{k}(m) - t)\phi'(t) dt = \int_{\frac{k}{m}}^{S_{k}(m)}(S_{k}(m) - t)\phi'(t) dt + \int_{S_{k}(m)}^{\frac{k+1}{m}}(S_{k}(m) - t)\phi'(t) dt.
\end{equation}
Then, considering the constant sign of $(S_{k}(m)-t)$ on $\biggl[\D \frac{k}{m}, S_{k}(m)\biggr]$, and on $\biggl[\D S_{k}(m),\frac{k+1}{m}\biggr]$, equation (\ref{II.1}) allows us to obtain:
$$
\|h\|^2m_2\int_{\frac{k}{m}}^{S_{k}(m)}\!\!(S_{k}(m) - t) dt \leq \!\int_{\frac{k}{m}}^{S_{k}(m)}\!\!(S_{k}(m) - t)\phi'(t) dt \leq \!\|h\|^2M_2\int_{\frac{k}{m}}^{S_{k}(m)}\!\!(S_{k}(m) - t) dt,
$$
and,
$$
\|h\|^2M_2\int_{S_{k}(m)}^{\frac{k+1}{m}}(S_{k}(m) - t) dt \leq \int_{S_{k}(m)}^{\frac{k+1}{m}}(S_{k}(m) - t)\phi'(t) dt \leq \|h\|^2m_2\int_{S_{k}(m)}^{\frac{k+1}{m}}(S_{k}(m) - t) dt.
$$
that lead to the next two inequalities:
\begin{equation}\label{FF9}
\int_{\frac{k}{m}}^{\frac{k+1}{m}}(S_{k}(m)-t)\phi'(t) dt \leq \|h\|^2M_2\int_{\frac{k}{m}}^{S_{k}(m)}(S_{k}(m)-t) dt + \|h\|^2m_2\int_{S_{k}(m)}^{\frac{k+1}{m}}(S_{k}(m)-t) dt,
\end{equation}
and,
\begin{equation}\label{FF10}
\int_{\frac{k}{m}}^{\frac{k+1}{m}}(S_{k}(m)-t)\phi'(t) dt \geq \|h\|^2m_2\int_{\frac{k}{m}}^{S_{k}(m)}(S_{k}(m)-t) dt + \|h\|^2M_2\int_{S_{k}(m)}^{\frac{k+1}{m}}(S_{k}(m)-t) dt.
\end{equation}
Since we also have the two following results:
$$
\D \int_{\frac{k}{m}}^{S_{k}(m)}(S_{k}(m)-t) dt = \frac{\lambda^2}{2m^2} \hs \mbox{ and } \int_{S_{k}(m)}^{\frac{k+1}{m}}(S_{k}(m)-t) dt = -\frac{(\lambda-1)^2}{2m^2}, \hs \mbox{where } \lambda \equiv mS_{k}(m)-k,
$$
inequalities (\ref{FF9}) and (\ref{FF10}) lead to:
\begin{equation}\label{FF11}
\frac{\|h\|^2}{2m^2}P_1(\lambda) \leq \int_{\frac{k}{m}}^{\frac{k+1}{m}}(S_{k}(m)-t)\phi'(t) dt \leq  \frac{\|h\|^2}{2m^2}P_2(\lambda),
\end{equation}
where two polynomials $P_1(\lambda)$ and $P_2(\lambda)$ are defined by:
$$
P_1(\lambda) \equiv m_2\lambda^2-(\lambda -1)^2M_2 \hs \mbox{ and } \hs P_2(\lambda) \equiv M_2\lambda^2-(\lambda -1)^2m_2.
$$
Keeping in mind that we want to minimize $\epsilon_{a,m+1}(h)$, we find that the value of $\lambda$ which minimizes the polynomial
$P(\lambda) \equiv P_2(\lambda)-P_1(\lambda)=(M_2-m_2)(2\lambda^2-2\lambda+1)$ is $\D\lambda = \frac{1}{2}$.\sa
Then, for this value of $\lambda$, (\ref{FF11}) becomes:
\begin{equation}\label{FF12}
\D \frac{\|h\|^2}{8m^2}(m_2-M_2) \leq \int_{\frac{k}{m}}^{\frac{k+1}{m}}(S_{k}(m)-t)\phi'(t) dt \leq  \frac{\|h\|^2}{8m^2}(M_2-m_2)\,.
\end{equation}
Finally, by summing over $k$ between $0$ to $m$, we have from (\ref{FF5}) and (\ref{FF12}):
\begin{equation}\label{epsilon1_final}
\frac{\|h\|}{8m}(m_2-M_2) \leqslant \epsilon_{a,m+1}(h) \leqslant \frac{\|h\|}{8m}(M_2-m_2).
\end{equation}
Using first the definitions of $\lambda$ and of $S_{k}(m)$, and  using also that the weights $\omega_{k}(m), (k=0,m),$
satisfy (\ref{poids_norme}), we have
\begin{equation}\label{Skm}
\forall m \in \mathbb{N}^{*}, \forall k \in [0,m[ : S_{k}(m) = \sum \limits_{j = 0}^{k} \omega_{j}(m) = \frac{1}{2m} + \frac{k}{m},
\end{equation}
and the corresponding weights $\omega_{k}(m)$ are equal to:
\begin{equation}\label{poids_sous_contraintes}
\omega_{0}(m) = \omega_{m}(m) = \frac{1}{2m}, \mbox{ and, } \hs\omega_{k}(m) = \frac{1}{m},\, (k=0,m-1)\,.
\end{equation}
This completes the proof of Theorem \ref{theorem_1}.
\end{prooff}
\noindent To illustrate this refined first-order expansion formula, let us derive (\ref{Generalized_Taylor}) when $m=2$, (that is to say with three points). In that case, we readily get:
$$
\D f(a+h) = f(a) + \left(\frac{D\!f(a) + 2D\!f\biggl(\D a+\frac{h}{2}\biggr)+ D\!f(a+h)}{4}\right)\!\!.(h) + \|h\|\epsilon_{a,3}(h),
$$
where 
\[ \frac{\|h\|}{16}(m_2-M_2) \leqslant \epsilon_{a,3}(h) \leqslant \frac{\|h\|}{16}(M_2-m_2).\]
\begin{remark}\label{Relation_Poids}
Condition (\ref{poids_norme}) on the weights $\omega_k(m), 0 \leq k \leq m,$ in Theorem \ref{theorem_1} is a kind of {\em closure condition} that helped us to determine $w_k(m)$. But it is not a restrictive one. Indeed, without this condition, one would have to consider (\ref{FF4}) in the place of (\ref{FF5}). \sa
Hence, $f$ being twice differentiable, it also exists $(m_1,M_1) \in\R^2$ such that
$$
\forall \, x\in U, \forall h\in \R^n: m_1\|h\|\hs \leq \hs D\!f(x).(h) \leq M_1\|h\|.
$$
Then, using that $\phi(1)= D\!f(a+h).(h)$, we obtain:
$$
m_1 \|h\| \leq \phi(1) \leq M_1 \|h\|,
$$
that leads to, together with (\ref{FF4}):
\begin{equation}\label{epsilon1_final_V2}
\frac{\|h\|}{8m}(m_2-M_2) - \frac{M_1}{2m} \leqslant \epsilon_{a,m+1}(h) \leqslant \frac{\|h\|}{8m}(M_2-m_2) - \frac{m_1}{2m}\,.
\end{equation}
Here, we used that the weights $\omega_k(m), (k=0,m),$ may be determined by the help of (\ref{Skm}) without taking into account anymore the closure condition (\ref{poids_norme}). \sa
More precisely, in that case, we get that the weights $\omega_k(m), (k=0,m),$ are equal to:
$$
\omega_{0}(m) = \frac{1}{2m} \hs \mbox{ and } \hs \omega_{k}(m) = \frac{1}{m}, \hs (k=1,m).
$$
Consequently, we deduce from (\ref{epsilon1_final_V2}) that the bounds of the reminder $\epsilon_{a,m+1}(h)$ are $m$ times lower than the ones given using the classical first-order Taylor's formula, see (\ref{epsilon1_bounded}).\sa
Finally, by considering the closure condition (\ref{poids_norme}) and the corresponding weights $\omega_k(m), (k=0,m)$, determined by (\ref{poids_sous_contraintes}),
we improved the result of (\ref{epsilon1_final_V2}), since the bounds of the remainder given by (\ref{FF4}) are $2m$ smaller than the ones given by the first Taylor's formula.
\end{remark}

\section{Application to the interpolation error}\label{applInter}

\subsection{The case of dimension one}
\noindent In this subsection, we consider the refined first-order expansion formula (\ref{Generalized_Taylor}) with two points, (i.e. $m=1$) in the one-dimensional case, namely with $U=]x_0,x_1[, (x_0<x_1)$. More precisely, formulas (\ref{Generalized_Taylor})-(\ref{majoration_epsilon_m+1}) give in this case, for $x_0 < a <x_1$:
\begin{equation}\label{Formule_2_points_V2.0}
\D f(a+h) = f(a) + \left(\frac{Df(a) + Df(a+h)}{2}\right)\!\!.(h) + \|h\|\epsilon_{a,2}(h),
\end{equation}
that can be simply written, in the one dimensional case, by using the derivative of $f$ and the absolute value as the norm in $\R$:
\begin{equation}\label{Formule_2_points_V2}
\D f(a+h) = f(a) + \left(\frac{f'(a) + f'(a+h)}{2}\right)\!\!.h + |h|\epsilon_{a,2}(h),
\end{equation}
with
$$
\frac{|h|}{8}(m_2-M_2) \leqslant \epsilon_{a,2}(h) \leqslant \frac{|h|}{8}(M_2-m_2).
$$
Now, for $x_0 < a <b  <x_1$, let $\Pi_{[a,b]}(f)$ be the usual interpolation polynomial of degree less than or equal to one defined by:
$$
\forall x \in [a,b]: \Pi_{[a,b]}(f)(x) = \frac{x-b}{a-b} \hs f(a) + \frac{x-a}{b-a} \hs f(b).
$$
Our aim is to investigate the consequences of formula (\ref{Formule_2_points_V2}) when we use it to estimate the error of interpolation $e(.)$ defined by
$$
\forall x \in [a,b]: e(x) = \Pi_{[a,b]}(f)(x) - f(x),
$$
and to compare it with the error obtained by the first order Taylor formula (\ref{Taylor_1_V1}),
written for $n=1$.
\begin{lemma}\label{error_estimate_1D}
Let $f$ be a function of $C^2(]x_0,x_1[)$. The following interpolation error estimate holds:
\begin{equation}\label{error_P1_M1_M2}
\forall x \in [a,b]: |\Pi_{[a,b]}(f)(x) - f(x)| \leq \frac{(b-a)}{4}\|f'\|_{\infty} + \frac{(b-a)^2}{16}\|f''\|_{\infty},
\end{equation}
where $\|.\|_{\infty}$ denotes the classical norm on $L^{\infty}$. \sa
\end{lemma}
\begin{prooff}
We showed in \cite{ChaskJams} (see Formula (56)) that the usual interpolation polynomial $\Pi_{[a,b]}(f)$ introduced above can be written, using  (\ref{Formule_2_points_V2}), as:
\begin{equation}\label{Error_P1_V00}
\D \Pi_{[a,b]}(f)(x) = f(x) + \biggl(\frac{f'(b)-f'(a)}{2(b-a)}\biggr)(b-x)(x-a) + \frac{(b-x)(x-a)}{(b-a)}\left[\frac{}{}\!\epsilon_{x,2}(b)-\epsilon_{x,2}(a)\right],
\end{equation}
where $\epsilon_{x,2}(b)$ and $\epsilon_{x,2}(a)$ are the remainders of (\ref{Formule_2_points_V2}) for $a$ and $b$ respectively.\sa
Consider now these two remainders written in the integral form. Due to (\ref{FF5}), for $m=1$, we have 
$$
\D \epsilon_{x,2}(a) =  \int_{0}^{1}\!\left(\frac{1}{2}-t\right)\!(a-x) f''\!\bigl(x+t(a-x)\bigr)\,dt,
$$
and
$$
\D \epsilon_{x,2}(b) =  \int_{0}^{1}\!\left(\frac{1}{2}-t\right)\!(b-x) f''\!\bigl(x+t(b-x)\bigr)\,dt.
$$
As a consequence,  using that $\D\int_{0}^{1}\left|\frac{1}{2}-t\right|dt = \frac{1}{4}$, these two remainders are bounded by
$$
\D |\epsilon_{x,2}(a)| \leq \frac{(x-a)}{4}\|f''\|_{\infty} \, \mbox{ and } \, |\epsilon_{x,2}(b)| \leq \frac{(b-x)}{4}\|f''\|_{\infty}\,.
$$
Now, using that $\D \sup_{a\leq x\leq b}(x-a)(b-x)=\frac{(b-a)^2}{4}$, together with (\ref{Error_P1_V00}) yields
\begin{equation}
\D|\Pi_{[a,b]}(f)(x) - f(x)| \leq \frac{(b-a)\|f'\|_{\infty}}{4} + \frac{(b-a)^2\|f''\|_{\infty}}{16}\,.
\end{equation}
\end{prooff}
Considering now the usual Taylor's formula, classically, the interpolation error $e(.)$ is bounded by (see for example \cite{Crouzeix_Mignot}, \cite{Ciarlet2})
$$
\D \forall x \in [a,b]: |f(x)-\Pi_{[a,b]}(f)(x)| \leq \frac{1}{2 }(x-a)(x-b)\|f''\|_{\infty}.
$$
So, using again that $\D \sup_{a\leq x\leq b}(x-a)(b-x)=\frac{(b-a)^2}{4}$, we get from this error bound that
\begin{equation}\label{error_interpolation_litterature_V2}
\D \forall x \in [a,b]: |f(x)-\Pi_{[a,b]}(f)(x)|  \leq \frac{(b-a)^2}{8}\|f''\|_{\infty}.
\end{equation}
Now, we aim at evaluating the improvement obtained in the upper bound involved in (\ref{error_P1_M1_M2}), compared to the one deduced from the usual interpolation error estimate (\ref{error_interpolation_litterature_V2}). Hence, we are looking for functions $f$ and for a positif number $\beta<1$ such that
\begin{equation}\label{Improve_2.0}
\frac{(b-a)\|f'\|_{\infty}}{4} + \frac{(b-a)^2\|f''\|_{\infty}}{16} \leq \beta\,\frac{(b-a)^2\|f''\|_{\infty}}{8},
\end{equation}
that can be rewritten as
\begin{equation}\label{Improve_2}
\D \|f'\|_{\infty} \leq \frac{\|f''\|_{\infty}}{\Lambda}, \quad \mbox{where } \Lambda = \frac{4}{(2\beta -1)(b-a)}\,.
\end{equation}
Hence, our goal is to identify values of $\beta$ and class of functions $f$ such that (\ref{Improve_2}) is satisfied. Obviously, (\ref{Improve_2}) does not hold for $0 \leq \beta \leq 1/2$.  In addition, since we are interested by illustrating (\ref{Improve_2}), we will also restrict ourselves by looking for functions $f$ solutions to:
\begin{equation}\label{Classe1}
\D \forall x \in [a,b]: |f'(x)| \leq \frac{f''(x)}{\Lambda},
\end{equation}
which implies that, in the sequel, we will only consider convex functions, i.e. functions such that $f''(x) \geq 0, \forall x \in[a,b]$, since it is clear that solutions of (\ref{Classe1}) also satisfy (\ref{Improve_2}).\sa
The next lemma determines a necessary condition for inequality (\ref{Classe1}) to be verified:
\begin{lemma}\label{lemma-illustrate-1d}
Let $f$ be a convex function of $C^2([a,b])$, which satisfies inequality (\ref{Classe1}). Then, we have
$$
\D \forall x \in [a,b], f(x) \geq \max\biggl(f(a) + \frac{f'(a)}{\Lambda}\bigl(e^{\Lambda(x-a)}-1\bigr), f(a) - \frac{f'(a)}{\Lambda}\bigl(e^{-\Lambda(x-a)}-1\bigr)\biggr)\,.
$$
\end{lemma}
\begin{prooff}
From differential inequality (\ref{Classe1}), we get the two following second order differential inequalities, for all $x \in [a,b]$:
\begin{equation*}
	\left \{ \begin{array}{lll}
		 f''(x) - \Lambda f'(x) & \geq &  \hs 0,  \\[0.2cm]
		 f''(x) + \Lambda f'(x) &  \geq  & \hs 0, 
	\end{array} \right.
\end{equation*}
that can be rewritten, setting $F(x)=f'(x)$, 
\begin{equation*}
	\left \{ \begin{array}{lll}
		F'(x) - \Lambda F(x) & \geq &  \hs 0,  \\[0.2cm]
		F'(x) + \Lambda F(x) &  \geq  & \hs 0, 
	\end{array} \right.
\end{equation*}
Now, multiplying the first inequality by the integrating factor $e^{-\Lambda(x-a)}$ and the second one by the integrating factor 
$e^{\Lambda(x-a)}$, we obtain that
$$
\D \frac{d}{dx}\biggl[F(x)e^{-\Lambda(x-a)}\biggr] \geq 0\,, \hs \quad \mbox{ and } \quad \hs \frac{d}{dx}\biggl[F(x)e^{\Lambda(x-a)}\biggr] \geq 0\,.
$$
This means that both functions $F(x)e^{-\Lambda(x-a)}$ and $F(x)e^{\Lambda(x-a)}$ are increasing on $[a,b]$, or in other words, returning to function $f'$:
$$
\D \forall x \in [a,b]: f'(x) \geq e^{\Lambda(x-a)}f'(a) \mbox{ and } f'(x) \geq e^{-\Lambda(x-a)}f'(a).
$$
It suffices to integrate thes two differential inequalities on the interval $[a,x], (x \in [a,b])$ to obtain the inequality of Lemma \ref{lemma-illustrate-1d}.
\end{prooff}
The next lemma enables us to determine sufficient conditions to determine a class of functions $f$ satisfying condition (\ref{Classe1}), and consequently condition (\ref{Improve_2.0}).
\begin{lemma}\label{Id_Class_Funct}
Let $\delta > 0$ be a real number and let $f$ be the solution to the problem \textbf{(P}\textbf{)} defined by, for all $x \in [a,b]$:
\begin{subequations}
\begin{empheq}[left=\mbox{\textbf{(P)}}\hs\empheqlbrace]{alignat=2}
\hs  f''(x)  & -  \Lambda f'(x)  =  \hs \delta,  &  & \label{EquaDiff_0} \\[0.2cm]
\hs  f'(a) &  \geq -\frac{\delta}{\Lambda}. &  &  \label{Cond_init}
\end{empheq}
\end{subequations}
Then, we have
\begin{equation}\label{Groupe_f}
\forall \beta\in\mbox{$]\frac{1}{2},1]$}: f(x) = f(a) + \frac{f'(a)}{\Lambda}\bigl(e^{\Lambda(x-a)}-1\bigr) + \frac{\delta}{\Lambda}\biggl[\frac{e^{\Lambda(x-a)}-1}{\Lambda}-(x-a)\biggr].
\end{equation}
In particular, $f$ given by (\ref{Groupe_f}) is a convex function which satisfies (\ref{Classe1}).
\end{lemma}
\begin{prooff}
Two successive integrations lead to the solution $f$ given by (\ref{Groupe_f}). Let us check now that these functions $f$ are convex. Their second derivative is equal to
$$
\D f''(x)=\biggl(\!f'(a)+\frac{\delta}{\Lambda}\biggr)\Lambda e^{\Lambda(x-a)}
$$
which is positive, since $f$ solution to problem \textbf{(P)} satisfies condition (\ref{Cond_init}).\sa
Furthermore, remark also that, due to condition (\ref{Cond_init}), $f$ solution to problem \textbf{(P)} satisfies as well
\begin{equation}\label{Ineq+Lambda}
f''(x) + \Lambda f'(x) \geq 0,  \forall x \in [a,b].
\end{equation}
Indeed, considering that $f$ is given by (\ref{Groupe_f}), we have, for all  $x \in [a,b]$:
$$
f''(x) + \Lambda f'(x) = 2\biggl(\!f'(a)+\frac{\delta}{\Lambda}\biggr)\Lambda e^{\Lambda(x-a)}-\delta.
$$
Hence, for all  $x \in [a,b]$, $f''(x) + \Lambda f'(x)\geq 0$ is equivalent to
$$\D e^{\Lambda(x-a)}\geq \frac{\delta}{2\biggl(f'(a)+\displaystyle\frac{\delta}{\Lambda}\biggr)\Lambda}, \forall x \in [a,b],$$
which leads to $\D f'(a)\geq -\frac{\delta}{\Lambda}$, using the minimum of the exponential for $x=a$.\sa
Finally, we obtained that solutions $f$ of (\ref{EquaDiff_0})-(\ref{Cond_init}) determined by (\ref{Groupe_f}) fulfill inequality (\ref{Ineq+Lambda}). In other words, solutions $f$ given by (\ref{Groupe_f}) satisfy inequality
(\ref{Classe1}).
\end{prooff}
\noindent With Lemma \ref{Id_Class_Funct}, we have checked that condition (\ref{Classe1}) does not lead to an empty set of functions. Indeed, functions determined by (\ref{Groupe_f}) satisfy condition (\ref{Classe1}), and then (\ref{Improve_2.0}), for any value of $\beta\in]\frac{1}{2},1]$. For these functions, the interpolation error $e(.)$ is $\beta$ times smaller than the one found with the usual first-order Taylor's formula. So, the best we can obtain by the refined first-order expansion formula (\ref{Formule_2_points_V2}) corresponds to a decrease of about 50 percents $(\beta=\frac{1}{2})$.

\subsection{The case of dimension $n,(n> 1)$}\label{Dimension_n}
\noindent Let us now consider formula (\ref{Formule_2_points_V2.0}) together with the related classical one (\ref{Taylor_1_V1}) in $\R^n, (n>1)$. \sa
First of all, let us rewrite both formulas by using the integral form of remainder. From propositions \ref{Prop1} and \ref{Prop2}, for any function 
$v \in C^2(U)$, the classical Taylor formula can be written as:
\begin{equation}\label{Taylor_Integral}
\D v(a+h) = v(a) + Dv(a).(h) + \int_{0}^{1}(1-t)\phi'(t)\,dt,
\end{equation}
where the function $\phi$ is defined in (\ref{phi}).\sa
On the other hand, let us consider formula (\ref{FF5}) of the remainder, corresponding to the refined first-order expansion formula with two points, that is for $m=1$. Using  (\ref{Skm}) and (\ref{poids_sous_contraintes}), we get the following integral form of (\ref{Formule_2_points_V2.0}):
\begin{equation}\label{Taylor_2pts_Integral}
\D v(a+h) = v(a) + \left(\frac{Dv(a) + Dv(a+h)}{2}\right)\!\!.(h) + \int_{0}^{1}\left(\frac{1}{2}-t\right)\phi'(t)\,dt.
\end{equation}
\noindent Let us now consider an open-bounded and non empty subset $\Omega$ of $\R^{n}$ which is simply connected, and let us denote by $\partial\Omega$ its boundary,
assumed to be a simplicial complex (i.e.  the generalization to $\R^n$ of a polygon in $\R^2$; \cite{Wallace}, Chapter IX).\sa
We  consider a generalized triangulation $\mathcal{T}_{\h}$ of $\bar{\Omega}$ composed by a finite number of $n$-simplicies $\mathcal{S}_k, (1 \leq k \le N),$ which respects the classical rules of a "\emph{regular}" discretization:
the set $\bar{\Omega}$ is expressed as the set-theoretic union of a finite number of $n$-simplices $\mathcal{S}_k, (1 \leq k \le N)$, whose interior are pairwise
disjoint, and such that, given any $n$-simplex of the triangulation, each one of its $(n-1)$-face is either a portion of the boundary $\partial\Omega$,
or an $(n-1)$-face of another $n$-simplex of the triangulation, (for more details, see for example \cite{Coatmelec}). \sa
For all $k\in[\![1,N]\!]$\footnote{as usual, the notation $[\![1,N]\!]$ denotes all the integers form $1$ to $N$}, we denote by $\h_k\equiv diam(\mathcal{S}_k)$ the diameter of the $n-$simplex $\mathcal{S}_k$, (the diameter being the greatest distance between two points inside
$\mathcal{S}_k$), and by $\h$ the mesh size ($\D\h=\max_{k=1,N}\h_k$) of the corresponding triangulation $\mathcal{T}_{\h}$.\sa
Finally, we introduce the vertices $A^{(k)}_i, i \in [\![1,n+1]\!],$ associated to a given $n-$simplex $\mathcal{S}_k, k \in [\![1,N]$. Then, we define by $\pi_{\h}(v)$ the
piecewise polynomial of degree less than or equal to one which belongs to $C^0(\bar{\Omega})$ such that
\begin{equation}\label{Pi_h_Def}
\forall k\in [\![1,N]\!] , \forall P\in S_k: \pi_{\h}(v)(P) = \pi_{\mathcal{S}_k}(v)(P),
\end{equation}
where, $\forall P \in S_k$, we have 
\begin{equation}\label{Pi_Sk}
\pi_{\mathcal{S}_k}(v)(P) = \sum_{i=1}^{n+1}\lambda_i(P)v(A^{(k)}_i)\,.
\end{equation}
Above, $\lambda_i$ are the barycentric functions which satisfy, $\forall P \in \mathcal{S}_k$, 
\begin{equation}\label{barycentric_functions}
\sum_{i=1}^{n+1}\lambda_i(P)OA^{(k)}_i \hs = \hs OP \hs \hs , \hs \hs \quad \mbox{ with } \sum_{i=1}^{n+1}\lambda_i(P) \hs = \hs 1,
\end{equation}
$O$ denoting the origin of $\R^n$.\sa
Therefore, $\pi_{\mathcal{S}_k}(v)$ is the unique polynomial of degree less than or equal to one which interpolates the function $v$ at the vertices $A_k, k \in [\![1,n+1]\!]$, (see \cite{CiarletWAG} for example).\sa
\noindent On a given $n-$simplex $\mathcal{S}$ (we drop the subscript $_{k}$ for simplicity), our aim is now to evaluate the distance between a "smooth function" $v$ and its interpolated polynomial $\pi_{\mathcal{S}}(v)$, first by the classical Taylor formula (\ref{Taylor_Integral}), then, by our new refined first-order expansion formula (\ref{Taylor_2pts_Integral}). \sa
To this end, let us recall some classical notations:
If $|||.|||$ denotes the norm of the operator $D^2(v)(P)$, $(\forall P \in \bar{S}),$ defined in (\ref{norme_forme_bili}) and the corresponding one
for $D(v)(P)$, and on the other hand, we define the $L^\infty$-norm $\|.\|_{\infty}$ for any integer $p$ and a $p$-linear mapping fields
$L: P \in \bar{S}\rightarrow {\cal{L}}\bigl((\R^n)^p;\R\bigr)$ by:
$$
\|L\|_{\infty} = \sup_{P\in \bar{S}}|||L(P)|||.
$$
Then, we have obtained the following interpolation error estimates:
\begin{lemma}\label{Thm_Interpolation_complete}
Let $\mathcal{S}$ be a given $n-$simplex in $\R^n$ defined by its vertices $A_i,i \in [\![1,n+1]\!]$. Let $v$ be a function of $C^2(\mathcal{S})$, we have the two following interpolation error estimates:
\begin{eqnarray}
\forall P \in \mathcal{S} & : &  |\pi_{\mathcal{S}}(v)(P)-v(P)| \, \leq \, \frac{\|D^2(v)\|_{\infty}}{2}diam(S)^2, \label{local_interpolation_error_1_V0} \\[0.1cm]
\forall P \in \mathcal{S} & : &  |\pi_{\mathcal{S}}(v)(P)-v(P)| \, \leq \, \frac{\|D(v)\|_{\infty}}{2}diam(S) + \frac{\|D^2(v)\|_{\infty}}{4}diam(S)^2\,.
\label{local_interpolation_error_2_V0}
\end{eqnarray}
\end{lemma}
\begin{prooff}
Let us begin to prove the estimate (\ref{local_interpolation_error_1_V0}) related to the classical Taylor's formula. Using formula (\ref{Taylor_Integral}) by taking for $a$ any point $P \in S$,  we obtain that
$$
v(A_i) = v(P) + D\!v(P).(h_i) + \int_{0}^{1}(1-t)\phi'_i(t)\,dt,
$$%
where the functions $\{\phi_i\}_{1 \leq i \leq n+1}$ are defined by:
$$
\forall t \in [0,1], \, \phi_i(t)=Dv(P+th_i).(h_i), \, \mbox{ and } \, h_i = PA_i.
$$
Let us compute now the quantity $\pi_{\mathcal{S}}(v)(P)$ defined by (\ref{Pi_Sk}). We have
\begin{eqnarray}
\pi_{\mathcal{S}}(v)(P) \!\! & \!\! = \!\! & \!\! \D \sum_{i=1}^{n+1}v(A_i)\lambda_i(P), \nonumber\\[0.2cm]
\!\! & \!\! = \!\! & \!\! v(P)\sum_{i=1}^{n+1}\lambda_i(P) + \sum_{i=1}^{n+1}\lambda_i(P)D\!v(P).(h_i) + \sum_{i=1}^{n+1}\lambda_i(P)\int_{0}^{1}(1-t)\phi'_i(t)\,dt, \nonumber\\[0.2cm]
\!\! & \!\! = \!\! & \!\! v(P) + D\!v(P).\!\!\left(\sum_{i=1}^{n+1}\lambda_i(P)h_i\!\!\right) + \sum_{i=1}^{n+1}\lambda_i(P)\int_{0}^{1}(1-t)\phi'_i(t)\,dt. \label{EQ1}
\end{eqnarray}
Above, we used that $D\!v(P)$ is a linear form together with the second property of (\ref{barycentric_functions}). Moreover, due to the two properties of (\ref{barycentric_functions}), we also have:
$$
\D\sum_{i=1}^{n+1}\lambda_i(P)h_i = \sum_{i=1}^{n+1}\lambda_i(P)OA_i-\sum_{i=1}^{n+1}\lambda_i(P)OP = OP - OP = 0.
$$
Then, (\ref{EQ1}) leads to:
$$
\D \pi_{\mathcal{S}}(v)(P) = v(P) + \sum_{i=1}^{n+1}\lambda_i(P)\int_{0}^{1}(1-t)\phi'_i(t)\,dt,
$$
that implies the following inequality, $\forall P \in S$, and using that, $\forall i \in [\![1,n+1]\!], 0 \leq \lambda_i \leq 1$:
\begin{equation}\label{EQ4}
 |\pi_{\mathcal{S}}(v)(P)-v(P)| \leq \frac{1}{2} \max_{i=1,n+1}\|\phi'_i\|_{\infty}\,.
\end{equation}
However, we also have the following property: $\forall t \in [0,1], \forall P \in  S, \forall i\in[\![1,n+1]\!]$,
\begin{eqnarray}
\D |\phi'_{i}(t)| & = & |D^2(v)(P+t.h_{i}).(h_{i},h_{i})| \hs \leq \hs |||D^2(v)(P+t.h_{i})|||\, diam(S)^2, \nonumber\\[0.2cm]
& \leq & \|D^2(v)\|_{\infty} \, diam(S)^2, \label{Diam2}
\end{eqnarray}
where we set above $h_{i}\equiv PA_i$.\sa
Then, using (\ref{Diam2}), inequality (\ref{EQ4}) leads to the following result, $\forall P \in S$: 
\begin{equation}
|\pi_{\mathcal{S}}(v)(P)-v(P)| \leq \frac{\|D^2(v)\|_{\infty}}{2}diam(S)^2\,,
\end{equation}
that proves the estimate (\ref{local_interpolation_error_1_V0}) of the lemma.\\

\noindent Let us prove now the estimate (\ref{local_interpolation_error_2_V0}) related to the refined first-order expansion formula. In a similar way, we consider now formula (\ref{Taylor_2pts_Integral}) to compute the interpolation polynomial $\PiS$ defined in (\ref{Pi_Sk}) at point $P$. We obtain that:
\begin{eqnarray}
\PiSP \hs = \hs v(P) & + & \D\frac{1}{2}\sum_{i=1}^{n+1}\lambda_i(P)D\!v(P).(h_i) + \frac{1}{2}\sum_{i=1}^{n+1}\lambda_i(P)D\!v(A _i).(h_i) \nonumber \\[0.2cm]
& + & \D\sum_{i=1}^{n+1}\lambda_i(P)\int_{0}^{1}\left(\frac{1}{2}-t\right)\phi'_i(t)\,dt,  \nonumber\\[0.2cm]
\hs = \hs v(P) & + & \frac{1}{2}\sum_{i=1}^{n+1}\lambda_i(P)\phi_i(0) + \sum_{i=1}^{n+1}\lambda_i(P)\int_{0}^{1}\left(\frac{1}{2}-t\right)\phi'_i(t)\,dt, \label{EQ2}
\end{eqnarray}
and finally, that 
\begin{equation}\label{EQ5}
\D|\pi_{\mathcal{S}}(v)(P) - v(P)| \leq \frac{1}{2}\D \max_{i=1,n+1}\|\phi_i\|_{\infty} + \frac{1}{4}\D \max_{i=1,n+1}\|\phi'_i\|_{\infty}.
\end{equation}
Moreover, we also have, for the function $\phi_{i}(t)$:
\begin{eqnarray*}
\D |\phi_{i}(t)| & = & |Dv(P+t.h_{i}).(h_{i})| \hs \leq \hs |||D(v)(P+t.h_{i})|||\, diam(S), \nonumber \\[0.2cm]
& \leq & \|D(v)\|_{\infty} \, diam(S),
\end{eqnarray*}
that yields, using inequality (\ref{EQ5}):
\begin{equation}\label{EQ6}
\D|\pi_{\mathcal{S}}(v)(P) - v(P)| \leq \frac{\|D(v)\|_{\infty}}{2}\D diam(S) + \frac{\|D^2(v)\|_{\infty}}{4}\D (diam(S))^2.
\end{equation}
\end{prooff}
To conclude this section, we can summarize the lemma (\ref{Thm_Interpolation_complete}) in the following way:  $\forall P \in \mathcal{S}$, the interpolation error is lower than the minimum between the upper bounds of (\ref{local_interpolation_error_1_V0}) and (\ref{local_interpolation_error_2_V0}), namely
\begin{equation}\label{New_Error_Estimate}
\D |\pi_{\mathcal{S}}(v)(P)-v(P)| \leq \min \left(\frac{\|D^2(v)\|_{\infty}}{2}diam(S)^2 , \frac{\|D(v)\|_{\infty}}{2}diam(S) + \frac{\|D^2(v)\|_{\infty}}{4}diam(S)^2\right).
\end{equation}

\noindent The next section is devoted to investigate consequences of our main result, regarding finite elements applications.

\section{Application to finite elements error estimates}\label{appliFE}
\noindent In this section, we consider the refined first-order expansion formula (\ref{Taylor_2pts_Integral}), namely  the general formula  (\ref{Generalized_Taylor}) for $m=1$.  Then, we study the impact of this formula in the context of Lagrange finite elements error estimate.\sa
We first recall the mathematical framework. Let us consider again a non empty open-bounded and simply connected
subset $\Omega$ of $\R^{n}$,  as introduced in section \ref{Dimension_n}, together with the generalized triangulation $\mathcal{T}_{\h}$ of $\bar{\Omega}$. Classically, we denote by $\mu(\Omega)$ the measure of $\Omega$ and by $V$ a Hilbert space, endowed with a norm $ \left\|.\right\|_{V}$, made of functions defined on $\Omega$. \sa
Then, we consider a linear continuous form $l(\cdot)$ defined on~$V$, and a bilinear, continuous and $V-$elliptic form $a(\cdot,\cdot)$ defined on $V \times V$. Particularly,  $\exists \,(\alpha,C)\in \R_+^*\times\R_+^*$ such that
\begin{equation}\label{Continuity_Coercivity}
\forall v \in V, \alpha \left\|v\right\|^{2}_{V} \leq a(v,v) \leq C \left\|v\right\|^{2}_{V}.
\end{equation}
Now, let $u \in V$ be the unique solution to the second order elliptic variational formulation \textbf{(VP)} defined by:
\begin{equation}\label{VP}
\textbf{(VP}\textbf{)} \hspace{0.2cm} \left\{
\begin{array}{l}
\mbox{Find } u \in V \mbox{ solution to:} \\ [0.2cm]
a(u,v) = l(v), \forall v \in V, \\[0.2cm]
\end{array}
\right.
\end{equation}
and let also introduce the approximation $u_{\h}$ of $u$, solution to the approximate variational formulation \textbf{(VP)}$_{\h}$:
\begin{equation}\label{VP_h}
\textbf{(VP}\textbf{)}_{\h} \hspace{0.2cm} \left\{
\begin{array}{l}
\mbox{Find } u_{\h} \in   V_{\h} \mbox{ solution to:} \\ [0.2cm]
a(u_{\h},v_{\h}) = l(v_{\h}), \forall v_{\h} \in V_{\h}, \\[0.2cm]
\end{array}
\right.
\end{equation}
where $V_{\h}$ denotes a finite-dimensional subset of $V$.\sa
The first step to estimate the error between $u$ and  $u_{\h}$ is given by C\'ea's Lemma \cite{ChaskaPDE}:
\begin{lemma}\label{Lemma_CEA}
Let $u$ denote the solution to (\ref{VP}) and $u_{\h}$ the solution to (\ref{VP_h}). Then, the following inequality holds:
$$
\left\|u - u_{\h} \right\|_{V} \leq \frac{C}{\alpha} \inf_{v_{\h} \in V_{\h}} \left\|u - v_{\h} \right\|_{V},
$$
where the constant $C$ and $\alpha$ are respectively  the continuity constant and the ellipticity constant of the bilinear form $a(\cdot,\cdot)$ defined in (\ref{Continuity_Coercivity}).
\end{lemma}
From C\'ea's lemma, it results that, to estimate the approximation error $\left\|u - u_{\h} \right\|_{V}$, we have to choose an element $v^{*}_{\h}\in V_{\h}$ for which an estimate of $\left\|u - v^{*}_{\h} \right\|_{V}$ can be computed. A convenient well-known choice consists in choosing $v^{*}_{\h}$ as an interpolation polynomial of a given degree.\sa
\noindent In the sequel of this section, we will consider several choices of $v^{*}_{\h}$. We will study the consequences of the  refined first-order expansion formula compared to the classical Taylor formula.
\subsection{The case of $V^{(1)}_{\h}$}

\noindent In this subsection, we consider the case when the Hilbert space $V$ is the Sobolev space $H^1(\Omega)$. As a first choice for $v^{*}_{\h}$, we assume that the approximated variational space $V_{\h}$ is equal to the finite dimensional polynomial subspace $V^{(1)}_{\h} \subset H^1(\Omega)$, defined by:
\begin{eqnarray}
V^{(1)}_h & = & \left\{v^{(1)}_{h}:\Omega\rightarrow\R, \,v^{(1)}_{h}\!\in C^{0}({\bar{\Omega}}), \restriction{v^{(1)}_{h}}{S_k} \!\in P_1(S_k)\right\},
\end{eqnarray}
where $P_{1}(S_k)$ denotes the set of polynomials defined on a given $n-$simplex $S_k$ whose degree is less than or equal to $1$.\sa
In this case, we have the following error estimate for $P_1$ finite element method:
\begin{theorem}\label{Thm_error_estimate_V1}
Let $\Omega \subset \R^n$ be an open bounded and simply connected domain,  and let $\mathcal{T}_h$ be a regular finite element mesh of $\Omega$. Assume that  the exact solution $u$ to (\ref{VP}) belongs to $C^2(\bar{\Omega})$ and let $u^{(1)}_h\in V^{(1)}_{\h}$ be the corresponding approximate solution to (\ref{VP_h}). \sa
Then, we have the following error estimate:
\begin{equation}
\left\|u-u^{(1)}_{\h}\right\|_{L^2(\Omega)} \leq \frac{C}{\alpha}\left\|u-\pi_{\h}(u)\right\|_{L^2(\Omega)} \leq \D\frac{C}{\alpha}\min\left(\frac{\|D^2(u)\|_{\infty}}{2}\h^2, \frac{\|D(u)\|_{\infty}}{2}\h + \frac{\|D^2(u)\|_{\infty}}{4}\h^2\right) \sqrt{\mu(\Omega)}\,.
\end{equation}
\end{theorem}
\begin{prooff}
Due to C\'ea's lemma , we have
$$
\forall v_{\h}\in V^{(1)}_{\h}: \left\|u-u^{(1)}_{\h}\right\|_{L^2(\Omega)} \hs \leq \hs  \frac{C}{\alpha}\left\|u-v_{\h}\right\|_{L^2(\Omega)}.
$$
So, we choose as the particular element $v^{*}_{\h}$, the $P_1$ interpolation function $\pi_{\h}(u)$ defined by (\ref{Pi_h_Def}) which belongs to $V^{(1)}_{\h}$. This allows us to evaluate the $L^2-$ norm of the quantity $u-\pi_{\h}(u)$ as follows:
\begin{equation}\label{In_1_0}
\|u-\pi_{\h}(u)\|_{L^2(\Omega)}^{2} = \sum_{\mathcal{S}_k \in \mathcal{T}_{\h}}\int_{\mathcal{S}_k}|u-\pi_{\h}(u)|^2d\Omega =\sum_{\mathcal{S}_k \in \mathcal{T}_{\h}}\int_{\mathcal{S}_k}\left|\frac{}{}\!u_{|_{\mathcal{S}_k}}-\pi_{_{\mathcal{S}_k}}(u_{|_{\mathcal{S}_k}})\right|^2d\Omega,
\end{equation}
where we used the property that, $\forall \, \mathcal{S}_k \in \mathcal{T}_{\h}$, $(\pi_{\h}\,u)\bigl|_{{\mathcal{S}_k}} = \pi_{_{\mathcal{S}_k}}(u_{|_{\mathcal{S}_k}})$.
Moreover, due to the interpolation error (\ref{New_Error_Estimate}), we also have that
\begin{equation}\label{In_1_1}
\left|\frac{}{}\!u_{|_{\mathcal{S}_k}}-\pi_{_{\mathcal{S}_k}}(u_{|_{\mathcal{S}_k}})\right| \leq \min \left(\frac{\|D^2(u)\|_{\infty}}{2}\hs ,\hs
\frac{\|D(u)\|_{\infty}}{2} + \frac{\|D^2(u)\|_{\infty}}{4}\right).
\end{equation}
As a consequence, inequality (\ref{In_1_0}) becomes:
\begin{equation}\label{IN1_10}
\|u-\pi_{\h}(u)\|_{L^2(\Omega)}^{2} \leq \sum_{\mathcal{S}_k \in \mathcal{T}_{\h}}\!\! \min\biggl(\frac{\|D^2(u)\|_{\infty}}{2} \,diam(S_k)^2, \frac{\|D(u)\|_{\infty}}{2} \,diam(S_k) +\frac{\|D^2(u)\|_{\infty}}{4} \,diam(S_k)^2\biggr)^{\!2} \!\mu(\mathcal{S}_k)\,.
\end{equation}
where $\mu(\mathcal{S}_k)$ denotes the measure of $\mathcal{S}_k$. \sa
Finally, from (\ref{IN1_10}) inequality (\ref{In_1_0}) leads to the following error estimate:
$$
\left\|u-u^{(1)}_{\h}\right\|_{L^2(\Omega)} \hs \leq \hs  \frac{C}{\alpha}\left\|u-\pi_{\h}(u)\right\|_{L^2(\Omega)} \leq  \frac{C}{\alpha}\min\left(\frac{\|D^2(u)\|_{\infty}}{2}\h^2, \frac{\|D(u)\|_{\infty}}{2}\h + \frac{\|D^2(u)\|_{\infty}}{4}\h^2\right) \! \sqrt{\mu(\Omega)}\,.
$$
\end{prooff}
Also in that case, this error estimate allows us to get a noticeable improvement of the upper bound of the approximation error (around 50 percents smaller), 
as soon as the minimum involved in this inequality is equal to $\D\frac{\|D(u)\|_{\infty}}{2}\h + \frac{\|D^2(u)\|_{\infty}}{4}\h^2$. \sa
\subsection{The case of $V^{(2)}_{\h}$}
\label{case-of-Vh2}
\noindent We propose now a second example of element $v^{*}_{\h}$. We choose, as Hilbert space $V$, the the finite
dimensional polynomial subspace $V^{(2)}_{\h} \subset H^1(\Omega)$ defined by:
$$
V^{(2)}_{\h}  =  \left\{v^{(2)}_{\h}:\Omega\rightarrow\R, v^{(2)}_{\h}\!\in C^{0}({\bar{\Omega}}), \restriction{v^{(2)}_{\h}}{S_k} \!\in P_2(S_k) \right\},
$$
where $P_{2}(S_k)$ denotes the set of polynomials defined on a given $n-$simplex $S_k$ whose degree is less than or equal to $2$.\sa
Let $\mathcal{S}$ denote a given $n-$simplex in $\R^n$. We first formulate an interpolation error result for a polynomial function of degree less than or equal to two, but for a function $v$ which only belongs to $C^2(\mathcal{S})$. For our purpose here, we introduce the {\em corrected }interpolation polynomial of $\pi_{\mathcal{S}}(v)$ defined by:
$$
\forall P \in T : \pi^{*}_{\mathcal{S}}(v)(P) = \pi_{\mathcal{S}}(v)(P)-\frac{1}{2}\sum_{i=1}^{n+1}\lambda_i(P)Dv(A_i).(PA_i).
$$
\begin{corollary}\label{Thm_Interpolation_Pi*}
Let $\mathcal{S}$ be a given $n-$simplex in $\R^n$ defined by its vertices $A_i, (i=1,n+1),$ and let $v$ be a function of $C^2(\mathcal{S})$. Then, the following interpolation error estimate holds:
\begin{eqnarray}
\forall P \in \mathcal{S} & : & |\pi^{*}_{\mathcal{S}}(v)(P)-v(P)| \, \leq \, \frac{\|D^2(v)\|_{\infty}}{4} \,\, diam(\mathcal{S})^2,\vspace{-0.4cm}  \label{local_interpolation_error_2}
\end{eqnarray}
\end{corollary}
\begin{prooff}
This is an immediate consequence of (\ref{EQ2}).
\end{prooff}

\noindent Now, if we compare the interpolation error estimate (\ref{local_interpolation_error_2}) to the classical one (\ref{local_interpolation_error_1_V0}), this naturally leads us to derive a new upper bound for the  finite elements error in $V^{(2)}_{\h}$. We have obtained the following error result:
\begin{theorem}\label{Thm_error_estimate}
Let $\Omega \subset \R^n$ be an open bounded and simply connected domain,  and let $\mathcal{T}_h$ be a regular finite element mesh of $\Omega$. Assume that the exact solution $u$ to (\ref{VP}) belongs to $C^2(\bar{\Omega})$ and let $u^{(2)}_h\in V^{(2)}_{\h}$ be the corresponding approximate solution to (\ref{VP_h}). \sa
Then, we have the following error estimate:
\begin{equation}\label{estimation_error_2}
\left\|u-u^{(2)}_{\h}\right\|_{L^2(\Omega)} \leq \frac{C}{\alpha}\left\|u-\pi^{*}_{\h}(u)\right\|_{L^2(\Omega)} \leq \D\frac{C\|D^2(u)\|_{\infty}}{4\,\alpha} \h^2\, \sqrt{\mu(\Omega)}\,.
\end{equation}
\end{theorem}
\begin{prooff}
Again, by C\'ea's lemma , we have
\begin{equation}\label{estimation_error_3}
\forall v_{\h}\in V^{(2)}_{\h}: \left\|u-u^{(2)}_{\h}\right\|_{L^2(\Omega)} \hs \leq \hs  \frac{C}{\alpha}\left\|u-v_{\h}\right\|_{L^2(\Omega)},
\end{equation}
and we choose for $v_{\h}$ the interpolation function $\pi^{*}_{\h}(u)$ introduced above. So we evaluate the $L^2-$ norm
of the quantity $u-\pi^{*}_{\h}(u)$ as we did in Theorem \ref{Thm_error_estimate_V1}:
$$
\|u-\pi^{*}_{\h}(u)\|_{L^2(\Omega)}^{2} = \sum_{\mathcal{S}_k \in \mathcal{T}_{\h}}\int_{\mathcal{S}_k}|u-\pi^{*}_{\h}(u)|^2d\Omega = \sum_{\mathcal{S}_k \in \mathcal{T}_{\h}}\int_{\mathcal{S}_k}\left|\frac{}{}\!u_{|_{\mathcal{S}_k}}-\pi^{*}_{_{\mathcal{S}_k}}(u_{|_{\mathcal{S}_k}})\right|^2d\Omega\,,
$$
that becomes, due to the interpolation error (\ref{local_interpolation_error_2}):
\begin{eqnarray}\label{umuh2bis}
\|u-\pi^{*}_{\h}(u)\|_{L^2(\Omega)}^{2} \hs & \leq & \hs \frac{\|D^2(u)\|_{\infty}^2}{16} \sum_{\mathcal{S}_k \in \mathcal{T}_{\h}} diam(\mathcal{S}_k)^4\,\mu (\mathcal{S}_k), \\[0.1cm]
\hs & \leq & \hs \frac{\|D^2(u)\|_{\infty}^2}{16} \h^4 \,\mu(\Omega). 
\end{eqnarray}
Finally, using this last estimate,  inequality (\ref{estimation_error_3}) leads to:
$$
\left\|u-u^{(2)}_{\h}\right\|_{L^2(\Omega)} \hs \leq \hs  \frac{C}{\alpha}\left\|u-\pi_{\h}(u)\right\|_{L^2(\Omega)} \leq \frac{C\|D^2(u)\|_{\infty}}{4\alpha}\, \h^2 \,\sqrt{\mu(\Omega)}.
$$
\end{prooff}
\begin{remark}
$\frac{}{}$
\begin{enumerate}
\item Using $P_{2}$ finite element, that is when $V_{\h}=V^{(2)}_{\h}$, one can get  \cite{Ciarlet2} an error bound smaller than what we got in (\ref{estimation_error_2}). However, this result is obtained by assuming an additional regularity to the exact solution $u$, namely $u \in C^3(\bar{\Omega})$. In our case, only the $C^2(\bar{\Omega})$-regularity is requested, and the classical result can not be applied anymore.

However, even if $u_{h}$ belongs to $V^{(2)}_{\h}$, we can also choose here $v^{*}_{\h}=\pi_{\h}(u)$, that is, the interpolate function of degree less than or equal to one, see (\ref{Pi_h_Def}). By performing the same computations as those used to derive the above theorem, we obtain that 
$$
\left\|u-u^{(2)}_{\h}\right\|_{L^2(\Omega)} \leq \frac{C}{\alpha}\left\|u-\pi_{\h}(u)\right\|_{L^2(\Omega)} \leq \D\frac{C\|D^2(u)\|_{\infty}}{2\, \alpha}
\h^2 \,\sqrt{\mu(\Omega)}\,,
$$
which corresponds to an upper bound that is  two times greater than those derived in (\ref{estimation_error_2}).\\
\item A practical consequence of Theorem \ref{Thm_error_estimate} is the possibility of using a coarser mesh for a given accuracy. 
Indeed, assume that we want to ensure the approximation error $\left\|u-u^{(2)}_{\h}\right\|_{L^2(\Omega)}$ to be less than or equal to a given $\epsilon>0$, $\epsilon$ specifying the expected accuracy. Let $\h$ (respectively $\h^*$) denote the mesh size required to get this accuracy with $\pi_{\h}(u)$ (respectively with $\pi^{*}_{\h}(u)$). Following (\ref{umuh2bis}) and (\ref{estimation_error_2}), this requires
$$
\frac{C\|D^2(v)\|_{\infty}}{2\,\alpha} \h^2 \leq \epsilon, \quad \mbox{ and } \quad \D\frac{C\|D^2(v)\|_{\infty}}{4\,\alpha} \h^{*2} \leq \epsilon.
$$
Hence, $\h^*$ can be $\sqrt{2}$ times greater than $\h$ for a given accuracy. Roughly speaking, it means that the same accuracy can be reached by the two approaches, but with a size mesh with approximatively $0.7$ times fewer nodes in each direction. For instance in dimension three, this allows us to use, for a given accuracy, a mesh with about $0.7^{3}\simeq 0.34$ which corresponds to about $2/3$ of nodes less than in the standard method.
\end{enumerate}
\end{remark}

\section{Conclusions and perspectives}\label{D}
\noindent In this paper we derived a refined first-order expansion formula in $\R^n$ to minimize the unknown remainder which appears in the classical Taylor's formula.  For a given function, this new formula is composed by a linear combination of its first derivatives, computed at $m+1$ equally spaced points. We showed that the corresponding remainder can be minimized for a suitable choice of the weights involved in this linear combination. In particular, we proved that the new remainder is $4m$ smaller than the one which appears in the classical first Taylor's formula.\sa
Afterwards, we considered two important applications: the interpolation error and the finite elements error estimates. In both cases, we showed that we can gain a significant improvement of the error estimate upper bounds. For example, in the one-dimensional case, when $m=1$ (with two points involved in the refined formula), we showed that the upper bound of these errors is four times smaller than the usual ones estimated by the classical Taylor formula.\sa
Concerning the finite elements error estimates, for linear second elliptic PDE's, since the approximation error is bounded by the interpolation error, we proved that, by the help of the {\em corrected} interpolation polynomial introduced in subsection (\ref{case-of-Vh2}), we obtained for the interpolation error, a lower  upper bound than the usual one.\sa
Several other applications can also be concerned by this new refined first-order expansion formula. For example, the approximation error involved  in ODE's approximation where Taylor's formula is basically used to derive numerical schemes. \sa

\noindent \textbf{\underline{Homages}:} The authors want to warmly dedicate this research to pay homage to the memory of Professors Andr\'e Avez and G\'erard Tronel who largely promote the passion of research and teaching in mathematics of their students.
\end{document}